\newcommand{\bull}{\vrule height .9ex width .8ex depth -.1ex}
 \newcommand{\ppp}{\hfill $\bull$ }
 \documentclass[11pt]{article}
\usepackage{epsfig}
 \author{ M.-L. Labbi\thanks{
  Address: Department of Mathematics, College
 of Science, University of Bahrain, Isa Town 32038 Bahrain.
  E-mail: labbi@sci.uob.bh }}
   \title{Manifolds  with positive  second H. Weyl curvature invariant}
   \date{}
   \usepackage{amsmath}
   \usepackage{amscd}
\newtheorem{theorem}{Theorem}[section]
\newtheorem{corollary}[theorem]{Corollary}
\newtheorem{lemma}[theorem]{Lemma}

\begin{document}
   \maketitle
   \begin{abstract} The second H. Weyl curvature invariant of a Riemannian
   manifold, denoted $h_4$, is the
    second curvature invariant which appears in the well known tube formula
     of H. Weyl. It coincides with the Gauss-Bonnet integrand in dimension 4.
       A crucial property of $h_4$ is that it is nonnegative
        for Einstein manifolds, hence it provides a  geometric obstruction
        to the existence of Einstein metrics in dimensions $\geq 4$, independently from the
sign of the Einstein constant.
This  motivates  our study  of the positivity of this invariant.
Here in this paper we prove many constructions of metrics with positive 
second H. Weyl curvature invariant, generalizing similar well known results
for
the scalar curvature.
   \end{abstract}
   \par\bigskip\noindent
 {\bf  Mathematics Subject Classification (2000).} 53C21, 53B20.
   \par\medskip\noindent
   {\bf Keywords.}  H. Weyl curvature invariants, Einstein manifold, surgery.
   \par
 \section{Introduction and statement of the results}
 Let $(M,g)$ be a smooth Riemannian manifold of dimension $n\geq 4$. Let
 $R$, $cR$ and $c^2R$ denote respectively the Riemann curvature tensor,
 Ricci tensor and the scalar curvature of $(M,g)$. The {\bf {\sl second
 Hermann Weyl curvature invariant}}, which throughout this paper
  shall be written in abridged form as {\bf {shwci}}  and denoted by $h_4$,
  can be defined by
  \begin{equation*}
h_4=\|R\|^2-\|cR\|^2+{1\over 4}\|c^2R\|^2
\end{equation*}
\par
A crucial property of $h_4$ is that it is nonnegative for Einstein manifolds
(see section 3 below), and so it provides a new geometric obstruction to the
 existence
of Einstein metrics independently from the sign of the Einstein constant.
 In particular, the  manifolds which do not admit any
 metric with positive {\sl shwci} cannot admit any Einstein metric.\par
Recall that in dimensions greater than 4, we do not know  any topological
 restriction for a manifold to be Einstein. If one requires that the Einstein
 constant to be positive, then one has two geometric obstructions
 $cR>0$ and $c^2R>0$. \par\medskip\noindent
It would be then with a great benefit to have a classification of 
  manifolds with positive {\sl shwci}.\par\medskip
Here in this paper, we shall inaugurate the study of the positivity
 properties
 of this important invariant.\par
 In section 2, we introduce and study in general the Hermann Weyl curvature
  invariants which
appear in the tube formula. Many examples are included.\par
In section 3, we study separately the case of the second invariant.
 We prove that
 it is nonnegative for Einstein manifolds and nonpositive for conformally
 flat manifolds with zero scalar curvature. The limit cases are also discussed.
 \par Then we prove theorem A which states 
  that  positive (resp. nonnegative) $p$-curvature implies positive (resp.
  nonnegative) {\sl shwci}
 where $p=[(n+1)/2]$. In particular, positive (resp. nonnegative) sectional
  curvature implies
 positive (resp. nonnegative) {\sl shwci}. Also if $n\geq 8$, positive
 (resp. nonnegative) isotropic curvature implies
 positive (resp. nonnegative) {\sl shwci}.\par
 Then one can apply our previous constructions in the class of manifolds with positive
$p$-curvature (see \cite{Lab2,Lab3,Lab4}) to get many examples of metrics with positive shwci.\par
In section 4 we prove the following useful theorem. It generalizes a similar result for
the scalar curvature.:\par\medskip\noindent
{\bf Theorem B.} {\sl  Suppose  that the total space $M$ of a Riemannian submersion
 is compact
and the  fibers ( with the induced metric)
 are with positive {\sl shwci} then the manifold M admits a
 Riemannian metric with positive {\sl shwci}.}
\par\medskip\noindent
The section is then ended with two applications of this theorem.\par\medskip
In  section 5, we prove the following stability theorem in the class of compact manifolds with positive
{\sl shwci}.:
\par\medskip\noindent
{\bf Theorem C.} {\sl If a manifold $M$ is obtained from a compact manifold $X$ by surgery in codimension
$\geq 5$, and $X$ admits a metric of positive {\sl shwci}, then so does $M$.\par
In particular, the connected sum of two compact manifolds of dimensions $\geq 5$ and each one is with positive
{\sl shwci} admits a metric with positive {\sl shwci}. }
\par\medskip\noindent
Theorem C generalizes a celebrated theorem of Gromov-Lawson and Schoen-Yau for the scalar curvature.\par\smallskip\noindent
As a consequence of the previous theorem we prove that there are no restrictions on the fundamental group of a compact manifold of dimension $\geq 6$ to carry a metric with positive
{\sl shwci}.\par\medskip
Finally, let us mention that  it would be interesting to prove, like in the case of the scalar curvature, that every manifold with
nonnegative {\sl shwci} that is not identically zero  admits a metric with positive shwci. 
\section{The H. Weyl curvature invariants}
Let
 $\Lambda^{*}M=\bigoplus_{p\geq 0}\Lambda^{*p}M$ denote the ring of differential
  forms
 on $M$, where $M$ is as above.  Considering the tensor product over the ring of  smooth functions,
  we define
 ${\cal D}= \Lambda^{*}M\otimes \Lambda^{*}M=\bigoplus_{p,q\geq 0}
  {\cal D}^{p,q}$ where $  {\cal D}^{p,q}= \Lambda^{*p}M \otimes
   \Lambda^{*q}M$.
  It is  graded associative  ring and  called the ring of
   double forms on $M$. \par
   The ring of curvature structures on $M$ (\cite{Kulk}) is the ring ${\cal C}=\sum_{p\geq 0}
   {\cal C}^p$ where ${\cal C}^p$ denotes symmetric elements in
   ${\cal D}^{p,p}$. We denote by ${\cal C}_1$ (resp. ${\cal C}_2$, ${\cal C}_0$) the subring
of curvature structures satisfying the first (resp. the second, both the first and second) Bianchi identity. \par
The standard inner product and the Hodge star operator $*$ on $\Lambda^{*p}M$ can be extended in a
standard way to ${\cal D}$ and they satisfy the following properties, see \cite{Lab1} for the proof:
\begin{equation}\label{formula:a}g\omega=*c*\omega \end{equation}
for all $\omega\in {\cal D}$, where c denotes the contraction map.
Also  for all $\omega_1, \omega_2\in {\cal D}$, we have
 \begin{equation}\label{formula:b}<g\omega_1,\omega_2>=<\omega_1,c\omega_2>\end{equation}
 that is the contraction map is the formal adjoint of the multiplication map by the metric
 $g$. Furthermore, we have for all $\omega_1, \omega_2\in {\cal D}^{p,q}$
\begin{equation}\label{formula:c}<\omega_1,\omega_2>=*(\omega_1.*\omega_2)=*(*\omega_1.\omega_2)\end{equation}
and
\begin{equation}\label{formula:d}**=(-1)^{(p+q)(n-p-q)}Id\end{equation}
Where $Id$ is the identity map on ${\cal D}^{p,q}$.\par\medskip\noindent
Next, we define the H. Weyl curvature invariants: \par\medskip\noindent
{\bf Definition.}
{\sl The $2q$-Hermann Weyl curvature invariant, denoted $h_{2q}$, is
the complete contraction of the tensor $R^q$, precisely,
  $$h_{2q}={1\over (2q)!}c^{2q}R^q$$
where $R^q$ denotes the multiplication of $R$ with itself $q$-times in the ring ${\cal C}$.
}
\par\medskip
Remark that $h_2={1\over 2}c^2R$ is one half  the scalar curvature and if $n$ is even then $h_n$ is
(up to a constant) the Gauss-Bonnet integrand.
 \par\medskip\noindent
Note that in \cite{Lab1},  it is proved that
 \begin{equation}\label{formula:e}h_{2q}=
*{1\over (n-2q)!}g^{n-2q}R^q \end{equation}

\subsection{Examples}
\begin{enumerate}
\item Let $(M,g)$ be with constant sectional curvature $\lambda$, then
$$R={\lambda\over 2}g^2\qquad {\text and}\qquad R^q={\lambda^q\over 2^q}g^{2q}$$
And therefore $h_{2q}$ is constant and equals to
$$h_{2q}=*{1\over (n-2q)!}g^{n-2q}R^q=*{\lambda^q\over 2^q(n-2q)!}g^{n}
={\lambda^qn!\over 2^q(n-2q)!}$$
In particular,
\begin{equation} \label{formula:f}h_4={n(n-1)(n-2)(n-3)\over 4}\lambda^2\end{equation}
\par\bigskip\noindent
\item Let $(M,g)$ be a Riemannian product of two Riemannian manifolds
 $(M_1,g_1)$ and $(M_2,g_2)$. If we index by $i$ the invariants of the
 metric $g_i$ for $i=1,2$, then
 $$R=R_1+R_2\quad {\rm and}\quad R^q=(R_1+R_2)^q=
 \sum_{i=0}^qC_i^qR_1^iR_2^{q-i}$$
 consequently, a starightforward calculation shows that
 \begin{equation*}
 \begin{split}
 h_{2q}&={c^{2q}R^q\over (2q)!}=\sum_{i=0}^qC_i^q {c^{2q}\over (2q)!}
 (R_1^iR_2^{q-i})\\
 &=\sum_{i=0}^qC_i^q {c^{2i}R_1^i\over (2i)!}
 {c^{2q-2i}R_2^{q-i}\over (2q-2i)!}\\
 &=\sum_{i=0}^qC_i^q (h_{2i})_1 (h_{2q-2i})_2
\end{split} \end{equation*}
In particular,
\begin{equation}\label{formula:g}
h_4=(h_4)_1+{1\over 2}scal_1scal_2+(h_4)_2
\end{equation}
where $scal$ denotes the scalar curvature.
\par\bigskip\noindent
\item Let $(M,g)$  be a hypersurface of the Euclidean space.  If $B$ denotes the second
fundamental form at a given point, then the Gauss equation shows that
$$R={1\over 2}B^2\qquad {\text and}\quad R^q={1\over 2^q}B^{2q}$$
Consequently, if $\lambda_1\leq \lambda_2\leq ...\leq \lambda_n$ denote the eigenvalues
of $B$, then the eigenvalues of $R^q$ are ${(2q)!\over 2^q}\lambda_{i_1}
\lambda_{i_2}...\lambda_{i_{2q}}$ where $ i_1<...<i_{2q}$.
Consequently,
 $$h_{2q}=
 {(2q)!\over 2^q}\sum_{1\leq i_1<...<i_{2q} \leq n}\lambda_{i_1}
 ...\lambda_{i_{2q}}$$
 So they coincide, up to a constant, with the symmetric functions in the eigenvalues
 of $B$.\par\noindent
 \par\medskip\noindent
 \item  Let $(M,g)$  be a  conformally flat manifold. Then it is well known that
 at each point
 of $M$,
 the Riemann curvature tensor is determined by a symmetric bilinear form $h$,
 in the sens that $R=g.h$. Consequently, $R^q=g^qh^q$.\par
 Let $\{e_1,..., e_n\}$ be an orthonormal basis of eigenvectors of $h$ and
 $\lambda_1\leq \lambda_2\leq ...\leq \lambda_n$ denote the eigenvalues of $h$.
 \par
 Then it is not difficult to see that all the tensors $R^q$ are also
  diagonalizable by the $2q$-vectors $e_{\scriptstyle i_1}\wedge ...\wedge e_{\scriptstyle i_{2q}}$, $i_1<...<i_{2q}$.
 Their eigenvalues are of the form
$$R^q(e_1\wedge ...\wedge e_{2q},e_1\wedge ...\wedge e_{2q})=
 (q!)^2\sum_{1\leq i_1<...<i_q \leq 2q}\lambda_{i_1}
 ...\lambda_{i_q}$$
 Consequently we get
\begin{equation*}
h_{2q}={(n-q)!q!\over (n-2q)!}
\sum_{1\leq i_1<...<i_q \leq n}\lambda_{i_1}
 ...\lambda_{i_q}
 \end{equation*}
\item Let $g_t=tg$ for $t>0$. If we index by $t$ the invariants of $g_t$ then
$$R_t=tR\qquad {\rm and} \qquad R_t^q=t^qR^q$$
and therefore
\begin{equation}
(h_{2q})_t={1\over t^q}h_{2q}
\end{equation}
\end{enumerate}
 \par\medskip\noindent
Let us now recall some other useful facts from \cite{Lab1} which shall be used later.\par
Following Kulkarni we call the elements in $\ker c\subset D^{p,q}$
effective elements of $D^{p,q}$, and shall be denoted by $E^{p,q}$.
\par\noindent
Recall  the following  orthogonal decomposition
of $D^{p,q}$:
\begin{equation}\label{formula:h}
D^{p,q}=E^{p,q}\oplus gE^{p-1,q-1}\oplus g^2E^{p-2,q-2}\oplus ...
\oplus g^rE^{p-r,q-r}
\end{equation}
where $r={\min \{p,q\}}$.\par\noindent
With respect to the previous decomposition, if $\omega=\sum_{i=0}^p
g^{p-i}\omega_i\in {\cal C}_1^p$ and $n=2p$, then (see \cite{Lab1})
\begin{equation}\label{formula:i}
*\omega=\sum_{i=0}^p  (-1)^i
g^{p-i}\omega_i
\end{equation}
Also let us recall the following lemma from \cite{Lab1}:
\begin{lemma}\label{lemma:gpq} Let $\omega_1\in E_1^r,\omega_2\in E_1^s$ be effectives then
\begin{equation*}
\begin{split}
<g^p\omega_1,g^q\omega_2>=& 0\quad {\text if}\quad (p\not =q)
\quad {\rm or}
\quad ( p=q \quad {\rm and}\quad r\not=s)\\
<g^p\omega_1,g^p\omega_2>=& p!\bigl(\prod_{i=0}^{p-1}(n-2r-i)\bigr)<\omega_1,\omega_2>
\quad {\text if}\quad p\geq 1
\quad {\text and}
\quad r=s
\end{split}
\end{equation*}
\end{lemma}
\section{The second H. Weyl curvature invariant}
With respect to the previous orthogonal decomposition \ref{formula:h},
 the Riemann curvature tensor decomposes to
  $R=\omega_2+g\omega_1+g^2\omega_0$,  where
\begin{equation*}
\begin{split}
\omega_0=&{1\over 2n(n-1)}c^2R\\
 \omega_1=&{1\over n-2}(cR-{1\over n}gc^2R)\\
  \end{split}
  \end{equation*}
  and $\omega_2$ is the Weyl tensor, it is defined by the previous decomposition
  of $R$.\par\smallskip
Corollary 6.5 in \cite{Lab1} shows that
\begin{equation}\label{formula:j}
h_4=
{1\over (n-4)!}[n!||\omega_0||^2-(n-2)!||\omega_1||^2+(n-4)!||\omega_2||^2]
\end{equation}
using lemma \ref{lemma:gpq} we can easily check that
\begin{equation}
\begin{split}
\|\omega_2\|^2&=\|R\|^2-{1\over n-2}\|cR\|^2+{1\over 2(n-1)(n-2)}\|c^2R\|^2\\
\|\omega_1\|^2&={1\over (n-2)^2}(\|cR\|^2-{1\over n}\|c^2R\|^2)\\
\|\omega_0\|^2&={1\over 4n^2(n-1)^2}\|c^2R\|^2\\
\end{split}
\end{equation}
and consequently using formula \ref{formula:j} we obtain another
 useful expression for $h_4$ as follows:
\begin{equation}\label{formula:k}
h_4=\|R\|^2-\|cR\|^2+{1\over 4}\|c^2R\|^2
\end{equation}
The folowing theorem was first proved  in \cite{Lab1}.
\begin{theorem}
Let $(M,g)$ be a Riemannian manifold of dimension $\geq 4$.
\begin{enumerate}
\item If $(M,g)$ is an Einstein manifold then $h_4\geq 0$. Furthermore $h_4\equiv 0$ if and only if
 $(M,g)$ is flat.
\item
If  $(M,g)$ is conformally flat with zero scalar curvature then
$h_4\leq 0$. Furthermore $h_4\equiv 0$ if and only if
 $(M,g)$ is flat.
\end{enumerate}
\end{theorem}
{\bf Proof.} If $(M,g)$ is conformally flat then $\omega_2=0$ and then
\begin{equation*}
\begin{split}
h_4=&{1\over (n-4)!}[n!||\omega_0||^2-(n-2)!||\omega_1||^2]\\
=&{n-3\over n-2}[{n\over 4(n-1)}||c^2R||^2-||cR||^2]
\end{split}
\end{equation*}
From which is clear that if $c^2R=0$ then $h_4\leq 0$ and $h_4\equiv 0$
if and only if the metric
is Ricci flat and hence is flat. This proves the first part of the
 theorem.\par
Next, if $(M,g)$ is Einstein  then $\omega_1=0$ and hence
\begin{equation*}\begin{split}
h_4=&{1\over (n-4)!}[n!||\omega_0||^2+(n-4)!||\omega_2||^2]\\
=&||R||^2+{n-4\over 4n}(c^2R)^2
\end{split}
\end{equation*}
From which it is clear that  $h_4\geq 0$ and $h_4\equiv 0$ if and only if the metric
is flat. This completes the proof of the theorem.
\ppp
\par\medskip\noindent
Recall that (see \cite{Lab3,Lab4}) the $p$-curvature of $(M,g)$, denoted $s_p$ for $1\leq p\leq n-2$, is a function defined on the $p$-Grassmanian bundle of the manifold. Its value at a
tangent $p$-plane $P$ is the avearge of the sectional curvatures of all 2-planes orthogonal
to $P$. In particular $s_0$ is the scalar curvature and $s_{n-2}$ is twice the sectional curvature.\par
The following theorem provides a relation between the positivity of the $p$-curvature
and the {\sl shwci}.:\par\medskip\noindent
{\bf Theorem A.} {\sl
Let $(M,g)$ be a Riemannian manifold of dimension $n\geq 4$ and with
nonnegative (resp. positive) $p$-curvature such that $p\geq {n\over 2}$,
 then the {\sl shwci} of $(M,g)$ is nonnegative (resp. positive).
 Furthermore, it
 vanishes if and only if the manifold is flat.}
\par\medskip\noindent
{\bf Proof.} Suppose $n=2(k+2)$ is even, $k\geq 0$.
Since
$$R=\omega_2+g\omega_1+g^2\omega_0$$
then
$$g^kR=g^k\omega_2+g^{k+1}\omega_1+g^{k+2}\omega_0$$
and by formula \ref{formula:i} we have
$$*g^kR=g^k\omega_2-g^{k+1}\omega_1+g^{k+2}\omega_0$$
On the other hand since $s_{k+2}\geq0$, then both the tensors $g^kR$ and $*g^kR$ are with positive
sectional curvature, hence
$$[g^k\omega_2+g^{k+2}\omega_0](e_{i_1},...,e_{i_{k+2}},e_{i_1},...,e_{i_{k+2}})
\geq g^{k+1}\omega_1(e_{i_1},...,e_{i_{k+2}},e_{i_1},...,e_{i_{k+2}})$$
 and
$$[g^k\omega_2+g^{k+2}\omega_0](e_{i_1},...,e_{i_{k+2}},e_{i_1},...,e_{i_{k+2}})
\geq-g^{k+1}\omega_1(e_{i_1},...,e_{i_{k+2}},e_{i_1},...,e_{i_{k+2}})$$
for all orthonormal vectors $e_{i_1},...,e_{i_{k+2}}$,
and therefore
 $$[g^k\omega_2+g^{k+2}\omega_0](e_{i_1},...,e_{i_{k+2}},e_{i_1},...,e_{i_{k+2}})
\geq |g^{k+1}\omega_1(e_{i_1},...,e_{i_{k+2}},e_{i_1},...,e_{i_{k+2}})|$$
Consequently, using formulas \ref{formula:e} and \ref{formula:c},
 we get
$$h_4=*{1\over (n-4)!}g^{n-4}R^2=*{1\over (2k)!}(g^kR.g^kR)={1\over (2k)!}<g^kR,*g^kR>$$
and hence using lemma \ref{lemma:gpq} and considering an orthonormal basis
diagonalizing $cR$, we obtain
\begin{equation*}
\begin{split}
(2k)!h_4=&<g^k\omega_2+g^{k+2}\omega_0,g^k\omega_2+g^{k+2}\omega_0>-<g^{k+1}\omega_1,g^{k+1}\omega_1>
\\
\geq &\sum_{\scriptstyle i_1<...<i_{k+2}} \bigl[(g^k\omega_2+g^{k+2}\omega_0)(e_{i_1},...,e_{i_{k+2}},e_{i_1},...,e_{i_{k+2}})\bigr]^2
-  ||g^{k+1}\omega_1||^2  \\
\geq & \sum_{\scriptstyle i_1<...<i_{k+2}}
\bigl[g^{k+1}\omega_1(e_{i_1},...,e_{i_{k+2}},e_{i_1},...,e_{i_{k+2}})\bigr]^2-
||g^{k+1}\omega_1||^2=0
\end{split}
\end{equation*}
The same proof works for strict inequality. Also it is clear that
if $s_{k+2}\geq 0$ and  $h_4\equiv 0$ then $s_{k+2}\equiv 0$ so that  the
metric is flat.\par
To complete the proof, note that if the dimension of the manifold
$n=2p+1\geq 5$ is odd then one can consider the
 product $M\times S^1$. It is of even dimension $2(p+1)$ and  with nonnegative (resp. positive)
$(p+1)$-curvature therefore by formula
\ref{formula:g} we have
$$h_4(M)=h_4(M\times S^1)\geq 0 (\text{ resp. } >0)$$
this completes the proof of the theorem.\ppp
\begin{corollary}
\begin{enumerate}
\item A Riemannian manifold of dimension $\geq 4$ and with
nonnegative (resp.  positive) sectional curvature is with
nonnegative (resp.  positive) {\sl shwci}.
 Furthermore, $h_4\equiv 0$
  if and only if the metric is flat.
\item A Riemannian manifold of dimension $\geq 8$ and with
nonnegative (resp.  positive) isotropic curvature is with
 nonnegative (resp.  positive) {\sl shwci}. 
Furthermore, $h_4\equiv 0$
  if and only if the metric is flat.
\end{enumerate}
\end{corollary}
{\bf Proof.} Straightforward since positive sectional curvature implies positive $p$-curvature
and positive isotropic curvature implies the positivity of the $p$-curvature for all $p\leq n-4$,
see \cite{Lab2}.\ppp
\par\smallskip\noindent
{\bf Remarks.} 
\begin{enumerate}
\item  If the dimension of the manifold is even, say $n=2q$, Hopf conjecture states that if the sectional curvature is positive then so is the Gauss-Bonnet integrand that is $h_{2q}$. Then one can ask the more general question:\par\smallskip
{\sl Does positive sectional curvature implies positive $h_{2k}$, for all $2\leq 2k\leq n$?}\par\smallskip
This is now true for $k=1,2$ by the previous theorem, and it remains an open question for $k\geq 3$.
\item  Theorem A generalizes a
 result of  Thorpe \cite{Thorpe} for the dimension $n=4$.
\end{enumerate}
\par\medskip
It results from the previous corollary that Lie groups with a biinvariant metric and normal homogeneous
Riemannian manifolds are with nonnegative {\sl shwci}. Furthermore
 using our previous results on the $p$-curvature \cite{Lab3},
\cite{Lab4} and the above theorem we can easily prove the following
 corollaries:
\begin{corollary}
\begin{enumerate}
\item Let $G$ be a compact connected Lie group with rank $r$  such that $r<[{\dim G+1\over 2}]$ endowed
with a biinvariant metric $b$ then $(G,b)$ is with positive {\sl shwci}.\par
In particular, if $G$ is simple then it is with positive {\sl shwci}.
\item  If $G/H$ is a normal homogeneous Riemannian manifold such that the rank r of $G$ satisfies
$r<[{\dim(G/H)+1\over 2}]$ then it is with positive {\sl shwci}.
\end{enumerate}
\end{corollary}
\begin{corollary}
If a compact manifold $M$ admits a smooth action of a compact connected simple Lie group
with rank $r$ satisfying $r>[{\dim M+1\over 2}]$ then it admits a metric with positive {\sl shwci}.
\end{corollary}
\section{Proof of theorem B}
Let  $ (M,g) $ and  $ \left(B,\check  g \right) $ be two Riemannian manifolds,
and let $ \pi  :  (M ,g)\rightarrow  (B ,\check g)$ a Riemannian submersion.
 We define, for every  $t \in {\bf R}$, a new Riemannian metric $g_t$
on the manifold  $M$ by  multiplying the metric $g$ by $t^2$ in the vertical
 directions. Recall that $\forall m\in M$, we have a natural orthogonal
 decomposition of the tangent space at m
 $$T_mM={\cal V}_m\oplus {\cal H}_m$$
 where  ${\cal V}_m$ is the tangent to the fiber at m and ${\cal H}_m$
 is the horizontal space, so that
 \begin{eqnarray*}
 g_t \mid {\cal V}_m &=&t^2g\\
g_t \mid {\cal H}_m &=& g\\
g_t({\cal V}_m, {\cal H}_m)&=&0
\end{eqnarray*}
Note that in this case, $\pi :
(M,g_t)
\rightarrow (B,\check g)$ is still a Riemannian submersion with the same
horizontal and vertical distributions (see \cite{Besse}, \cite{Lab3}).\par \medskip
\noindent
In the following we shall index  by $t$ all the invariants of the metric
$g_t$, and in the case case
 $t=1$ we omitt the index 1.\par
Also, We make under a hat `` $\hat {} $ '' (resp. under a
 check `` $ \check {} $ ''
 ) the invariants of the fibers with the induced metric
 (resp. of the  basis  $B$ ).\par\medskip\noindent
Using lemma 2.1 in \cite{Lab3} it is easy to show that for all $g_t$-unit tangent vectors
$e_1,e_2,e_3,e_4$ we have
$$R_t(e_1,e_2,e_3,e_4)=O({1\over t})\qquad
\text {if one of these vecotors is horizontal}$$
and that
$$R_t(e_1,e_2,e_3,e_4)={1\over t^2}{\hat R}(te_1,te_2,te_3,te_4)+O(1)\qquad
{\text{ if the four vectors are vertical}}$$
Consequently, if $\{e_1,e_2,...,e_n\}$ is a $g_t$-orthonormal basis such that $\{e_1,...,e_q\}
\in {\cal V}_m$ and $\{e_{q+1},...,e_n\}
\in {\cal H}_m$, then
\begin{equation*}
\begin{split}
(||R_t||_t)^2=&\sum_{1\leq i<j\leq n,1\leq k<l\leq n}[R_t(e_i,e_j,e_k,e_l)]^2\\
=& {1\over t^4}\sum_{1\leq i<j\leq q,1\leq k<l\leq q}[{\hat R}(te_i,te_j,te_k,te_l)]^2+
O({1\over t^2})\\
=&{1\over t^4} \|\hat R \|^2+O({1\over t^2})\\
(||Ric_t||_t)^2=&\sum_{1\leq i,j\leq n}[Ric_t(e_i,e_j)]^2\\
=& {1\over t^4}\sum_{1\leq i,j\leq q}[{\hat Ric}(te_i,te_j)]^2={1\over t^4}||\hat Ric||^2+
O({1\over t^2})\\
(||scal_t||_t)^2=&{1\over t^4}||\hat scal||^2+O({1\over t^2})\\
\end{split}
\end{equation*}
Therefore, at the point $m$ we have:
\begin{equation}\label{h4MMM}
(h_4)_t= {1\over t^4}\hat{h_4}+O({1\over t^2})
\end{equation}
This completes the proof of theorem B since the total space is compact. \ppp
\par\medskip\noindent
 \begin{corollary}
\begin{enumerate}
\item The product $S^p\times M$ of an arbitrary compact manifold M
 with a sphere $S^p,p\geq 4$ admits a Riemannian metric with positive {\sl shwci}.
\item If a compact manifold admits a Riemannian foliation such that the
 leaves are with
  positive {\sl shwci} then the manifold admits a Riemannian metric with
  positive {\sl shwci}.
\end{enumerate}
\end{corollary}
\par\medskip\noindent
  {\bf Proof.} The first part  is straightforward, to prove the second one it suffices to
  notice that the proof of the previuos theorem works also in the case of local
  Riemannian submersions. \ppp
\par\medskip\noindent
 \begin{corollary}  If a compact manifold M admits a free and smooth action
  of a compact connected Lie group $G$ with rank $r$ such that $r<[{\dim G+1\over 2}]$
 then the manifold M admits a Riemannian metric with
  positive {\sl shwci}. \end{corollary}
\par\medskip\noindent
  {\bf Proof.} The canonical projection $M\rightarrow M/G$ is in this case a
  smooth submersion. Let the fibers be equipped with a biinvariant metric from
  the group $G$ via the canonical inclusion ${\cal G}\subset T_mM$.
  \par\noindent
  Using any G-invariant metric on M, we define the horizontal distribution
  to which we lift up an arbitrary metric from the basis $M/G$. Thus
  we have defined a metric on M such that the projection $M\rightarrow M/G$
  is a Riemannian submersion.\par\noindent
  Finally, since the group $G$ with a biinvariant metric is
  with positive {\sl shwci} then so are the fibers with the induced metric,
  and we conclude using the previous theorem. \ppp \par\medskip\noindent
{\bf Remark.} All simple Lie groups satisfy the property $r<[{\dim G+1\over 2}]$.

\section{Proof of Theorem C}
We proceed as in Gromov-Lawson's proof for the case of scalar curvature \cite{GroLaw}.\par
Let $(X,g)$ be a compact $n$-dimensional Riemannian manifold with positive {\sl shwci} and let
$S^m\subset X$ be an embedded sphere of codimension $q$ and with trivial normal bundle $N\equiv
S^m\times {\bf R}^q$. There exists $r_0>0$ such that the exponential map ${\rm exp}:
S^m\times D^q(r_{0})\rightarrow X$ is an embedding, where for every $x\in S^m$,  $\{x\} \times D^q(r_{0})$ denotes the closed Euclidean ball in ${\bf R}^q\equiv \{ x\}\times {\bf R}^q$. Let ${\rm exp}^*g$ denotes the pull back of
the metric $g$ to the normal sub-bundle  $S^m\times D^q(r_{0})$. 
\par\smallskip
Another natural metric on the normal bundle is the metric  $g^\nabla$
 defined using the normal connection $\nabla$, that is the metric compatible
  with the normal connection and such that the projection
  $ S^m\times D^q(r)\rightarrow S^m $
is a Riemannian submersion.
  We shall  denote also by $g^\nabla$ its restriction to the sub-bundles
 $S^m\times D^q(r)$ and  $\partial (S^m\times D^q(r))=S^m\times S^{q-1}(r)$.
 \par\smallskip
Recall that at each $(p,v)\in  S^m\times D^q(r)$ we have a natural
$g^\nabla$-orthogonal decomposition of the tangent space
into vertical and horizontal subspaces, namely
 \begin{equation}\label{riem:subm} T_{(p,v)} S^m\times D^q(r)
 ={\cal V}_{(p,v)}+{\cal H}_{(p,v)} \end{equation}
 where ${\cal V}_{(p,v)}$  is the tangent space
 to the fiber (over $p$) $D^q(r)$ at $v$.
These two metrics are tangent to the order two in the directions tangent to
 $D^q$, precisely
with respect to the decomposition \ref{riem:subm} we have (see \cite{Lab5})
\begin{equation}\label{bbb} 
\left( \begin{matrix} \displaystyle g^\nabla + 0(r^2 ) &\displaystyle
 & \displaystyle & g^\nabla + 0(r ) \cr\displaystyle
 g^\nabla + 0(r) &\displaystyle  & \displaystyle  &  g^\nabla + 0(r )\cr
\end{matrix}\right) 
\end{equation}
\noindent
{\bf Remark.}  Note that in \cite{GroLaw} in the begining of the proof of Lemma 2 page 430, it is claimed that
the former metrics are sufficiently close in the $C^2$-topology. But in general  this is only true for the directions
tangent to $S^{q-1}(r)$, a detailed study of the behavior of these two metrics will appear in a separate forthcoming  paper \cite{Lab5}. The same error is also in \cite{Lab4}. A  short proof of this fact is as follows:\par
With respect to the metric  $g^\nabla$, the sphere $S^m \hookrightarrow S^m\times D^q$ is totally geodesic (since for a Riemannian submersion the horizontal lift of a geodesic is a geodesic). But on the other side, the sphere  $S^m \hookrightarrow S^m\times D^q$ is totally geodesic
for the metric ${\rm exp}^*g$ only if the sphere $S^m$ is totally geodesic in $(X,g)$.\par\smallskip
However this does  not affect the corresponding conclusions in both papers
(after minor changes)
since the curvatures in question
 (that is the scalar curvature and the $p$-curvatures, $p\leq q-3$) of these
two metrics on the bundles $S^m\times S^{q-1}(r)$ are high and close enough
as $r\rightarrow 0$.\par\medskip
Now, it is easy to see that the second fundamental form of $S^m\times S^{q-1}(r)$
 in
$S^m\times D^q(r)$ with respect to the decomposition \ref{riem:subm}
 is of the form
\begin{equation}\label{aaa}
 \left( \begin{matrix} -{{\rm Id} \over r} &
    0 \cr   
0 &     0 \cr
\end{matrix} \right) 
\end{equation}

Consequently, using formulas \ref{bbb} and \ref{aaa} one can deduce without
 difficulties that
the second fundamental form of $S^m\times S^{q-1}(r)$ in
$S^m\times D^q(r)$ with respect to the  metric ${\rm exp}^*g$  is of the form (with respect to the decomposition \ref{riem:subm}):
\begin{equation}\label{ccc}
\left( \begin{matrix}  - {Id \over r}  + O(r) &
 &  O(1) \cr  &  &
 \cr O(1) &   &  O(1)\cr
 \end{matrix} \right) \end{equation}
\par
Note that since the second fundamental form is a continuous function, then it still has the form \ref{ccc} with respect to the following ${\rm exp}^{*}g$-orthogonal decomposition:
\begin{equation}\label{riem:submbis} T_{(p,v)} S^m\times D^q(r)
 ={\cal V}_{(p,v)}\oplus {\cal H}'_{(p,v)} \end{equation}
 where ${\cal V}_{(p,v)}$  is as in \ref{riem:subm} and the distribution ${\cal H}'$ is defined
by the previous orthogonal decomposition. Note that   as $r\rightarrow 0$, the distribution
${\cal H}'$ converges  to the distribution 
${\cal H}$ defined by the decomposition \ref{riem:subm}.
\par\medskip
Now we define a hypersurface $M$ in the  product
$ S^m\times D^q(r_{0})$ endowed with the product metric ${\rm exp}^{*}g \times {\bf R}$ by
 the relation
\[
M=\left\{ ((x,v),t)\in S^m\times D^q(r_{0})\times {\bf R}\quad /\quad (\Vert v\Vert ,t)\in
\gamma \right\}
\]
where $\gamma $ is a curve whose graph in the $(r,t)$-plane as pictured
below:
\begin{figure}[bht]
\begin{center}
\epsfig{figure=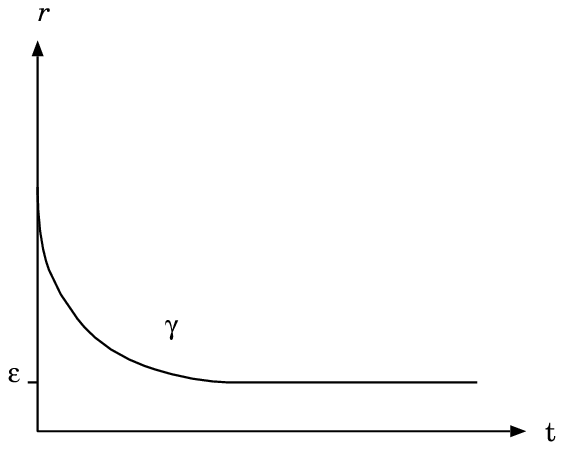,height=50mm,width=57mm,clip=,angle=0,
silent=,bbllx=6mm,bblly=245mm,bburx=63mm,bbury=293mm}
\end{center}
\end{figure}

The important points about $\gamma $ is that it is tangent to the $r$-axis
at $t=0$ and is constant for $r=\varepsilon >0$. Thus the induced metric on
$M$ extends the metric ${\rm exp}^{*}g$ on $S^m\times D^q(r_{0})$ near its boundary and
finishes with the product metric
$\bigl(\partial (S^m\times D^q(\varepsilon )),{\rm exp}^{*}g\bigr)
\times {\bf R}=\bigl(S^m\times  S^{q-1}(\varepsilon ),{\rm exp}^{*}g\bigr) \times
{\bf R} $.\par\medskip
Next, we evaluate the {\sl shwci} of the hypersurface $M$.\par
For each $m\in M$, we have the following ${\rm exp}^{*}g$-orthogonal decomposition of $T_mM$
\begin{equation}\label{ddd}
T_mM={\bf R}\tau \oplus {\cal V}_m\oplus {\cal H}'_m
\end{equation}
where $\tau$ is the unit tangent vector to the curve $\gamma$ in the $(r,t)$-plane and 
${\cal V}_m, {\cal H}'_m$ are as in \ref{riem:submbis}.\par\medskip
It results by a straightforward computation  using \ref{ccc} that the the second fundamental form of the hypersurface $M$
 has the following
 the form
(with respect to the decomposition \ref{ddd})
\begin{equation}\label{matrix}
\left( \begin{matrix}
k &  0 & ... &   0 \cr 
0 &  &  &    \cr
 &   &  \left(- {Id \over r}  +
 O(r)\right) {\sin} \theta  &   
 \bigl( O(1){\sin} \theta \bigr)\cr 
 \vdots & 
  &     &  \cr 
   & &  \bigl( O(1){\rm sin} \theta \bigr)
     & \bigl( O(1){\rm sin} \theta \bigr)\cr
    0 &  &    &
     \cr \end{matrix}\right)
    \end{equation}
where $k$ denotes the curvature of the curve $\gamma$ in the $(r,t)$-plane
and $\theta$ denotes the angle between the normal to $M$ and the $t$-axis at the corresponding point.\par\medskip
Then a long but direct computation using the Gauss equation and the previous formula \ref{matrix} shows that the curvatures of $M$ have the form :
\begin{equation*}
\begin{split}
\| R^M\|^2=&\|R^{\scriptstyle {S^p}\times {D^q}}\|^2+{(q-1)(q-2)\over 2r^4}\sin^4\theta+(q-1){k^2\over r^2}\sin^2\theta+O({1\over r^2})\sin\theta\\
\|{\rm Ric}^M\|^2=&\|{\rm Ric}^{\scriptstyle{ S^p\times D^q}}\|^2+{(q-1)(q-2)^2\over r^4}\sin^4\theta+q(q-1){k^2\over r^2}\sin^2\theta\\
&{\hfill -{(q-1)(q-2)^2k\over r^3}\sin^3\theta +O({1\over r^2})\sin\theta}\\
\|{\rm scal}^M\|^2=&\|{\rm scal}^{\scriptstyle S^p\times D^q}\|^2+{(q-1)^2(q-2)^2\over r^4}\sin^4\theta+4(q-1)^2{k^2\over r^2}\sin^2\theta\\
&{\hfill -2{(q-1)^2(q-2)k\over r^3}\sin^3\theta +O({1\over r^3})\sin\theta}\\
\end{split}
\end{equation*}
where we supposed that the curve $\gamma$ has its curvature of the form $k=O({1\over r})$.\par
\noindent
Consequently, we can evaluate the {\sl shwci} of $M$ as follows
\begin{equation}\label{h4M}\begin{split}
h_4^M=h_4^{\scriptstyle S^p\times D^q}+&{(q-1)(q-2)(q-3)(q-4)\over 4r^4}\sin^4\theta\\
&-{(q-1)(q-2)(q-3)k\over 2r^3}\sin^3\theta +O({1\over r^3})\sin\theta
\end{split}
\end{equation}
Next we shall show that it is possible to choose the curve $\gamma$ so that the metric induced
on $M$ has positive {\sl shwci} at all points $m\in M$.\par\smallskip
Formula \ref{h4M} shows that for $\theta=0$ we have $h_4^M=h_4^{\scriptstyle S^p\times D^q}$ is
positive, and then there exists an angle $\theta_0>0$ such that for all $0<\theta\leq \theta_0$
the {\sl shwci} of $M$ is positive.\par
then we continue with a straight line ($k=0$)  of angle $\theta_0$, say $\gamma_1$, until the term
${(q-1)(q-2)(q-3)(q-4)\over 4r^4}\sin^4\theta_0$ is strongly dominating.\par
On the other hand, when $\theta=\pi/2$ then $k=0$ and $r=\epsilon$ we have
\begin{equation}\label{h4MM}
h_4^M={(q-1)(q-2)(q-3)(q-4)\over 4\epsilon^4}+O({1\over \epsilon^3})
\end{equation}
which is positive as $\epsilon$ is small enough and $q\geq 5$.\par\smallskip

We now choose $r_1>0$ small and consider the point $(r_1,t_1)\in \gamma_1$. Then we bend
the straight line $\gamma_1$, begining at this point, with a curvature $k(s)$ of the following form
\begin{figure}[bht]
\begin{center}
\epsfig{figure=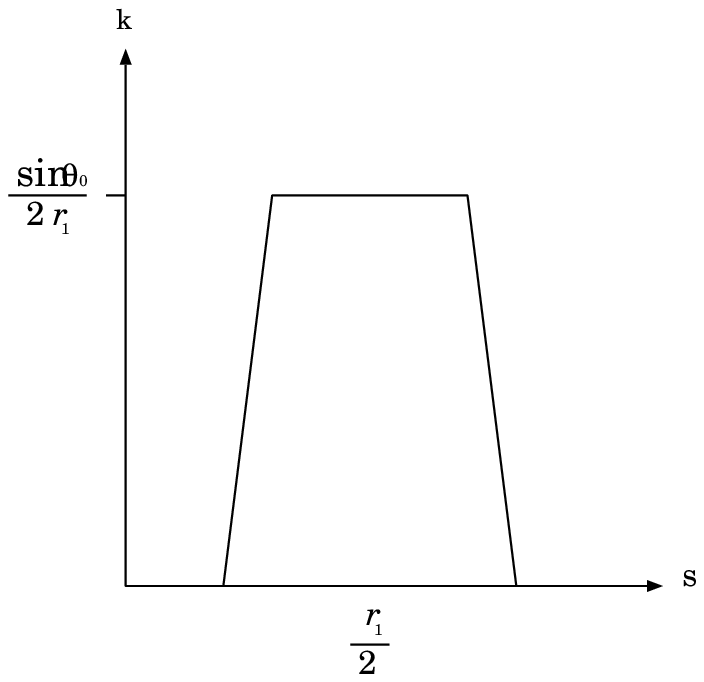,height=70mm,width=71mm,clip=,angle=0,
silent=,bbllx=2mm,bblly=222mm,bburx=73mm,bbury=293mm}
\end{center}
\end{figure}

where the variable $s$ denotes the arc length along the curve. \par
Since $q\geq 5$, formula \ref{h4M} shows that
\begin{equation}\label{eee}
h_4^M\geq h_4^{\scriptstyle S^p\times D^q}+{(q-1)(q-2)(q-3)\over 2r^3}\sin^3\theta
({\sin\theta\over 2r}-k)+O({1\over r^3})\sin\theta
\end{equation}
Then it is clear that 
that the hypersurface $M$ will continue to have $h_4^M>0$, since $k<\frac{\sin\theta_0}{2r_1}
<\frac{\sin\theta}{2r}$.\par
 After this first bending, we have
$\Delta r\leq \Delta s={r_1\over 2}$ and then $r\geq r_1-\Delta r\geq r_1-{r_1\over 2}>0$,
consequently the curve will not cross the $t$-axis.\par
On the other hand, $\Delta \theta =\int kds\approx {\sin\theta_0\over 4}$ is independent
of $r_1$. Clearly, by scaling down the curvature $k$, we can produce any
$\Delta \theta$ such that $0<\Delta \theta\leq  {\sin\theta_0\over 4}$.\par
Our curve now continues with a new straight line $\gamma_2$ with angle
 $\theta_1=\theta_0
+\Delta \theta$. By repeating this process finitely many times we can
 achieve a total bend
of ${\pi\over 2}$.\par\smallskip
Let $g_\epsilon$ denote the induced metric from ${\rm exp}^*g$ on
 $\partial(S^p\times D^q(\epsilon))=S^p\times  S^{q-1}(\epsilon) $,
and recall that the new metric defined  on $M$ is the old metric when $t=0$,
and finishes with the product
metric $g_\epsilon\times {\bf R}$. In the following we shall deform the
product metric $g_\epsilon \times {\bf R}$ on $ S^p\times  S^{q-1}(\epsilon) $,
 to the standard product metric through metrics with positive {\sl shwci}. This will
 be done in two steps:\par
Step1: We deform  The metric $g_\epsilon$ on
$S^m\times  S^{q-1}(\epsilon) $  to the standard product metric
$S^m(1)\times  S^{q-1}(\epsilon) $ through metrics with positive {\sl shwci}, as follows:\par
First, the metric $g_\epsilon$ can be homotoped through metrics
with $h_4>0$ to the normal metric $g^\nabla$ since their {\sl shwci} are
respectively high and close enough, see formulas \ref{h4MM} and \ref{h4MMM}.
\par
Then, for $\epsilon$ small enough, we can deform the normal metric $g^\nabla$
 on
$S^m\times  S^{q-1}(\epsilon )$ through Riemannian submersions to a new
metric where
$S^p$ is the standard sphere $S^p(1)$, keeping the horizontal distribution
 fixed.\par\noindent
This deformation keeps $h_4>0$ as far as $\epsilon$ is small enough,
 see formula \ref{h4MMM}.\par\smallskip
Finally, we deform the horizontal distribution to the standard one and
again by the same formula \ref{h4MMM}
this can be done keeping $h_4>0$. \par\medskip
Step2: Let us denote by $ds_t^2, 0\leq t\leq 1$, the previous family of
deformations on $S^m\times  S^{q-1}(\epsilon) $. They are all with positive
{\sl shwci}.  Where $ds_0=g_\epsilon$ and $ds_1$ is the standard product metric.\par
\noindent
It is clear that the metric $ds_{{t\over a}}^2+dt^2, 0\leq t\leq a$, glues together
the two  metrics $ds_0\times {\bf R}$ and $ds_1\times {\bf R}$. Furthermore,
there exists $a_0>0$ such that for all $a\geq a_0$ the metric
$ds_{{t\over a}}^2+dt^2$ on $S^m\times  S^{q-1}(\epsilon)\times [0,a] $
is with positive {\sl shwci}. In fact, via a change of variable, this is equivalent
to the existence of $\lambda_0>0$ such that for all $0<\lambda\leq \lambda_0$,
the metric $\lambda^2ds_t^2+dt^2$ is with positive {\sl shwci}. This is already known to be true
again by formula \ref{h4MMM}. \ppp
\begin{corollary} Let G be a finitely presented group.
Then for every $n\geq 6$, there exists a compact n-manifold M with positive
{\sl shwci} such that $\pi_1(M)=G$.
\end{corollary}
{\bf Proof.} Let $G$ be a group which has a presentation
consisting of $k$ generators $x_{1},x_{2},...,x_{k}$ and $l$ relations 
$r_{1},r_{2},...,r_{l}$.
\par\noindent Let $S^{1}\times S^{n-1}$ be endowed with the standard product
metric which is with positive {\sl shwci} ($n-1\geq 4$). Remark that the
fundamental group of $S^{1}\times S^{n-1}$ is infinite cyclic. Hence by
taking the connected sum $N$ of $k$-copies of $S^{1}\times S^{n-1}$ we
obtain an orientable compact $n$-manifold with positive {\sl shwci} (since
this operation is a surgery of codimension $n\geq 5$). By Van-Kampen
theorem, the fundamental group of $N$ is a free group on $n$-generators,
which we may denote by $x_{1},x_{2},...,x_{k}$.
\par\noindent We now perform surgery $l$-times on the manifold $N$,
killing in
succession the elements $r_{1},r_{2},...,r_{l}$. The result will be a
compact, orientable $n$-manifold $M$ with positive {\sl shwci} (since the
surgery is of codimension $n-1\geq 5$) such that $\pi _{1}(M)=G$, as
required. \ppp

\end{document}